\documentclass[11pt]{article}
\usepackage{amsfonts,amssymb}
\usepackage{amsmath}
\usepackage{amsthm}
\usepackage[matrix,arrow,curve]{xy}
\newcommand{\p}{\partial}
\newcommand{\e}{\varepsilon}

\newcommand{\R}{{\mathbb R}}
\newcommand{\E}{\overrightarrow{\exp}\int}
\newtheorem{theorem}{Theorem}
\newtheorem{lemma}{Lemma}
\newtheorem{prop}{Proposition}
\newtheorem{corollary}{Corollary}
\newtheorem{definition}{Definition}

\title{``Good Lie Brackets'' for Control Affine Systems}
\author{A.~A.~Agrachev\thanks{SISSA, Trieste}}
\date{}
\begin{document}
\maketitle

\begin{abstract}
We consider a smooth system of the form $\dot q=f_0(q)+\sum\limits_{i=1}^ku_if_i(q)$, $q\in M,\ u_i\in\R,$ and study
controllability issues on the group $\mathrm{Diff}M$. It is well-known that the system can arbitrarily well approximate
the movement in the direction of any Lie bracket polynomial of $f_1,\ldots,f_k$. Any Lie bracket polynomial of
$f_1,\ldots,f_k$ is good in this sense. Moreover, some combinations of Lie brackets which involve the drift term $f_0$
are also good but surely not all of them. In this paper we try to characterize good ones and, in particular, all universal good combinations, which are good for any nilpotent truncation of any system.
\end{abstract}

\section{Introduction}

We consider a control affine system
$$
\dot q=f_0(q)+\sum\limits_{i=1}^ku_if_i(q),\quad q\in M,\ u_i\in\R.  \eqno (1)
$$
Here $f_j,\ j=0,\ldots,k$, are smooth vector fields on the connected manifold $M$. In this paper, smooth means $C^\infty$. A smooth vector field $f$ is a first order differential operator on the space of smooth functions: $(fa)(q)=\langle d_qa,f(q)\rangle,\ \forall\,a\in C^\infty(M),\, q\in M$.
 The space of smooth vector fields forms a Lie algebra $\mathrm{Vec}M$ over $\R$, where Lie bracket $[f,g],\ f,g\in\mathrm{Vec}M$, is commutator of the differential operators, $[f,g]=f\circ g-g\circ f$.

Let $\mathrm{Lie}(f_0,\ldots,f_k)$ be Lie subalgebra of $\mathrm{Vec}M$ generated by $f_0,f_1,\ldots,f_k$. We assume that
$M$ is equipped with a complete Riemannian metric and all vector fields in $\mathrm{Lie}(f_0,\ldots,f_k)$ have at most linear growth in this metric; in particular, all these vector fields are complete.

Let $u(\cdot)=(u_1(\cdot),\ldots,u_k(\cdot))\in L_1([0,t];\mathbb R^k),\ f_u=\sum\limits_{i=1}^ku_if_i$. For any $q_0\in M$, there exists a unique solution
$q(\tau;u(\cdot)),\ 0\le\tau\le t$, of the equation
$
\dot q=f_0(q)+f_{u(\tau)}(q)
$
such that $q(0;u(\cdot))=q_0$. Moreover, the map
$$
q(0;u(\cdot))\mapsto q(t;u(\cdot)) \eqno (2)
$$
is a diffeomorphism of $M$. We address the controllability problem on the group of diffeomorphisms, it concerns the characterization of diffeomorphisms which can be realized, at least approximately, in this way.

\medskip\noindent{\sl Remark.}
We can treat control system as an evolutionary machine, a way to transform a linear time structure into the space structure.
Indeed, any control function $u(\cdot)$ (that is a function of time) produces a flow, a family of transformations of the space $M$.
Zero control provides a prescribed nominal dynamics (no events).
A change of control parameters means an event: the fields $f_0+u_if_i$ correspond to a short list of available simple immediate events.

The movement in the direction of a convex combination of the fields can be approximated by fast switchings among them:
coefficients in the convex combination are proportions of time used to follow the correspondent fields. Much more directions can be achieved by a clever use of the noncommutativity of the fields. It happens that extremely complex transformations can be approximated by applying available elementary events many times and in various orders.

The flows generated by the fields $f_0+u_if_i$ may be very simple and predictable; the complexity of the resulting transformation comes from the complexity of the evolutionary strategy, from the sophisticated choice of the order and timing for the switching between the fields based on the structure of their iterated Lie brackets.
This is why I say that control system is a machine to transform time structures into the space structures and why Lie bracket is a crucial tool.

\medskip

In what follows, we use chronological notations \footnote{See \cite[Ch.~6]{ABB} for a short self-contained presentation.
This chapter can be read independently on the rest of the book.} In particular, the diffeomorphism (2) is denoted by
$\E_0^tf_{u(\tau)}\,d\tau$. Attainable sets ${\mathcal A}_t\subset\mathrm{Diff}M$ are defined as follows:
$$
{\mathcal A}_t=\left\{\E_0^tf_0+f_{u(\tau)}\,d\tau: u(\cdot)\in L_1([0,t];\mathbb R^k)\right\}.
$$

A vector field $V\in\overline{\mathrm{Lie}(f_0,\ldots,f_k)}$ with at most linear grows is called a {\it good bracket} if
$e^{tV}\in\overline{\mathcal A}_t$ for any $t>0$. Here $\overline{S}$ is the closer of the set $S$; the closures of $\mathrm{Lie}(f_0,\ldots,f_k)$ and of ${\mathcal A}_t$
are taken in the $C^\infty$ topology. Good brackets provide us with additional control parameters, although these parameters are constrained. The following proposition is a simple corollary of the standard relaxation technique (see \cite[Sec.\,10.5]{AS}).

\begin{prop} Let $V$ be a good bracket, $u(\cdot)\in L_1([0,t];\mathbb R^k)$ and $v(\cdot)$ be a measurable real function on $[0,t]$ with values in $[0,1]$. Then the diffeomorphism
$$
\E_0^tv(\tau)(f_0+f_{u(\tau)})+(1-v(\tau))V\,d\tau
$$
belongs to $\overline{{\mathcal A}}_t.\qquad\square$
\end{prop}

In the next section, we construct an extensive family of good brackets. The construction is universal, it does not depend on the nature of vector fields $f_0,\ldots,f_k$. It actually demonstrates that certain control variations move the system
approximately in the direction of certain Lie bracket polynomials or series. The structure of these polynomials or series
depends only on the control variations while $f_0,\ldots,f_k$ serve as blind variables.

The whole construction can be performed in the free Lie algebra with free generators $a_0,\ldots,a_k$; this is done in section~3.
Theorem~4 in section~3 states that our construction describes all good brackets in the free Lie algebra. In other words,
the constructed family contains {\it all} Lie bracket series that are good brackets for {\it any} control affine system
such that the correspondent series of vector fields converge well.

Free Lie algebra is infinite dimensional and perhaps too big for practical applications. Finite dimensional truncations of the universal good brackets are studied in Section~4.

Any Lie algebra is quotient of a free one by an ideal. If projection of the set of universal good brackets to
the quotient algebra contains a nontrivial affine subspace, then the whole construction can be iterated and
more good brackets arise. This phenomenon is studied in Section~5.

In the last section~6, abstract results are applied to three special families of control systems.

Controllability paradigm is classical in the geometric control theory. More or less successful attempts to characterize various types of controllability in terms of iterated Lie brackets accompany the theory starting from its early days.
A part of this story can be recovered from \cite{J,S,K,AG}. The current paper is mainly inspired by recent works \cite{AS1,AS2}
on the ensemble controllability.

\medskip\noindent{\sl Acknowledgments.} I am grateful to Alexander Zuyev and Khazhgali Kozhasov for very useful stimulating discussions, and to Laboratoire Jacques-Louis Lions in Sorbonne Paris for the hospitality.

%We study controllability issues and do not impose constraints on the control parameters $u_i$ to simplify the task. We %however do not want to abuse this freedom too much and use only uniformly bounded trajectories to attain points of a %bounded domain.

\section{Fast-oscillating Control}

Let us start with small time control. To this end we take a sample control $u(t)=(u_1(t),\ldots,u_k(t)),\ t\in[0,2]$, a small parameter $\e>0$ and cook a re-scaled control $\frac 1\e u(\frac t\e)$. Now re-scaling time in the equation
$$
\dot q=f_0(q)+\frac 1\e f_{u(\frac t\e)}(q),
$$
we obtain the system
$$
\frac{dq}{d\tau}=\e f_0(q)+f_{u(\tau)}(q),
$$
where $\tau=\frac t\e$.

We have:
$$
q(2\e)=q_0\circ\overrightarrow{\exp}\int_0^2\e f_0+f_u\,d\tau=
$$
$$
 q_0\circ
 \overrightarrow{\exp}\int_0^2\e\overrightarrow{\exp}\int_0^\tau\mathrm{ad}f_{u(\theta)}\,d\theta f_0\,d\tau
 \circ\overrightarrow{\exp}\int_0^2f_u\,d\tau.
 $$

Now assume that $u(t)$ has a form:
$$
u(t)=\begin{cases}v(\tau),&0\le\tau\le 1,\\ -v(2-\tau),&1<\tau\le 2\,; \end{cases}
$$
then $\E_1^2 f_{u(\tau)}\,d\tau=\left(\E_0^1f_{u(\tau)}\,d\tau\right)^{-1}$ and we get:
$$
q(2\e)=q_0\circ\E_0^2\e\E_0^\tau\mathrm{ad}f_{u(\theta)}\,d\theta f_0\,d\tau=
$$
$$
q_0\circ e^{\e\left(\int_0^2\E_0^\tau\mathrm{ad}f_{u(\theta)}\,d\theta f_0\,d\tau+O(\e)\right)}.
$$

Moreover, $\E_0^{1+s}f_{u(\tau)}\,d\tau=\E_0^{1-s}f_{u(\tau)}\,d\tau,\ 0\le s\le 1$, and we obtain:
$$
\int\limits_0^2\E_0^t\mathrm{ad}f_{u(\tau)}\,d\tau f_0\,dt=2\int\limits_0^1\E_0^t\mathrm{ad}f_{v(\tau)}\,d\tau f_0\,dt,
$$
where $v(\cdot)$ is any function from $L_1([0,1];\R^k)$.
We set:
$$
{\cal V}=\left\{\int_0^1\E_0^t\mathrm{ad}f_{v(\tau)}\,d\tau f_0\,dt : v\in L_1([0,1];\R^k)\right\},
$$
substitute $2\e$ with $\e$ and we obtain:
$$
q(\e)=q_0\circ e^{\e\left(V+O(\e)\right)},\quad V\in\cal V.
$$

Let $t\mapsto V_t\in\cal V$ be a piecewise constant family of vector fields, $V_t=V_i,\ \forall t\in(\frac {i-1}\e,\frac i\e]$. We repeat our procedure $n$ times with an appropriate choice of $v(\cdot)$ for any segment $[\frac {i-1}\e,\frac i\e]$ and obtain:
$$
q(n\e)=q_0\circ e^{\e\left(V_\e+O(\e)\right)}\circ\cdots\circ e^{\e\left(V_{n\e}+O(\e)\right)}=
q_0\circ\E_0^{n\e}V_t+O(\e)\,dt.
$$

It follows that we can arbitrarily well approximate any diffeomorphism of the form $\E_0^tV_\tau\,d\tau$ by the value at the time moment $t$ of a flow generated by control system (1). Moreover, $t\mapsto V_t$ maybe of class $L_1$, it is not obliged to be piecewise constant.

Of course, we can also approximate any diffeomorphism of the form
$$
\E_0^tV_\tau\,d\tau\circ\E_0^tf_{u(t)}\,dt,\quad u(t)\in\R^k,\ V_\tau\in{\cal V}.
$$
The first term in the product is generated by the control linear system (without drift).

Assume for the moment that the fields $f_1,\ldots,f_k$ generate a finite-dimensional Lie algebra,
$\mathrm{Lie}\{f_1,\ldots,f_k\}=L,\ \dim L<\infty$. Let ${\cal L}\subset\mathrm{Diff}M$ be the Lie group generated by $L$;
it is a finite-dimensional Lie subgroup of $\mathrm{Diff}M$. Moreover,
$$
{\cal L}=\left\{\E_0^1f_{u(t)}\,dt : u(\cdot)\in L_1([0,1];\R^k)\right\}.
$$
.
\begin{prop} $\overline{{\cal V}}=\overline{\mathrm{conv}\{p_*f_0: p\in{\cal L}\}}.$
\end{prop}

\noindent {\bf Proof.} Recall that $p_*=\mathrm{Ad}p^{-1}$ and $\mathrm{Ad}\E_0^tf_\tau\,d\tau=\E_0^t\mathrm{ad}f_\tau\,d\tau$.

A convex combination of the vector fields $p_*f_0,\ p\in{\cal L},$ can be written as $\int_0^1\mathrm{Ad}p_tf_0\,dt$,
where $t\mapsto p_t$ is a piecewise constant family of diffeomorphisms from $\mathcal L$. We are interested in the closer of
the convexification and we may substitute piecewise constant families by smooth families or by any other class of families that is dense in $L_1([0,1];{\mathcal L})$.
The result now follows from the fact that any starting from $I$ curve in ${\mathcal L}$ can be uniformly approximated by the curve of the form $p_t=\E_0^tf_{u(\tau)}\,d\tau.\qquad \square$

In what follows, use a shorten notation:
$$
(\mathrm{Ad}{\mathcal L})f_0=\{(\mathrm{Ad}p)f_0: p\in{\mathcal L}\}=\{p_*f_0: p\in{\mathcal L}\}.
$$

\begin{corollary} Any diffeomorphism of the form
$$
\E_0^tV_\tau+X_\tau\,d\tau,\quad V_\tau\in\mathrm{conv}(\mathrm{Ad}{\mathcal L})f_0,\ X_\tau\in L,\ 0\le\tau\le t,
$$
belongs to $\overline{{\mathcal A}}_t.$
\end{corollary}

Indeed,
$$
\E_0^tV_\tau+X_\tau\,d\tau=\E_0^tW_\tau\,d\tau\circ p_t,
$$
where $p_\tau=\E_0^\tau X_\theta\,d\theta$ and $W_\tau=(\mathrm{Ad}p_\tau)V_\tau\in\mathrm{conv}(\mathrm{Ad}{\mathcal L})f_0$. Moreover, the curve $\tau\mapsto p_{\tau}$ in ${\mathcal L}$ can be uniformly approximated by a curve of the form $\tau\mapsto\E_0^tf_{u(\tau)}\,d\tau.\qquad\square$

\medbreak
We see that attainable sets of the system
$$
\dot q=U(q),\quad U\in\mathrm{conv}(\mathrm{Ad}{\mathcal L})f_0+L  \eqno (3)
$$
are contained in the closure of the attainable sets of system (1). Now if $\overline{\mathrm{conv}\{p_*f_0: p\in{\mathcal L}\}}$  contains an affine subspace and this affine subspace is not contained in $L$, then we can repeat the whole procedure and obtain more
available vector fields in the right hand side.

Of course, Proposition~2 can be extended to many interesting families of vector fields which generate infinite dimensional Lie algebras but we do not plan to do these functional analytic exercises in the current paper (cf.~\cite[Th.\,3]{AS1}). We are mainly interested in the structural questions.

Roughly speaking, the construction provides an extension of the original affine subspace of admissible vector fields to
a convex set in the closer of the Lie subalgebra generated by this affine subspace. The system can be moved in the direction of any field from this convex set that is built of Lie bracket polynomials and series of the original fields $f_j$. Moreover, the structure of the available polynomials and series does not depend on the choice of the original fields, the fields $f_j$ serve just as variables.

In other words, we may speak about a universal convex set in the closure of the free Lie algebra with generators
$a_0,a_1,\ldots,a_k$. This is the set of ``good'' combinations of brackets and any control affine system can approximate the movement in the direction of these combinations.

\section{Universal System}

The goal of this paper is to give as effective as possible description of ``good'' bracket combinations. To this end, we consider the universal control affine system in the closure of the universal associative algebra and finite-dimensional truncations of this system.

Let $\mathrm{Ass}(a_0,\ldots,a_k)$ be the free associated algebra over $\R$, its elements are linear combinations of words in the alphabet $\{a_0,\ldots,a_k\}$. Then \linebreak $\mathrm{Ass}(a_0,\ldots,a_k)=\bigoplus\limits_{n=0}^\infty A_n$, where the space $A_n$ consists of linear combinations of words with $n$ letters and $A_0=\R$ corresponds to the empty word.
The {\it closure} of $\mathrm{Ass}(a_0,\ldots,a_k)$ is the algebra of formal series
$$
{\frak A}=\left\{\sum\limits_{n=0}^\infty x_n: x_n\in A_n\right\}
$$
endowed with topology of the term-wise convergence; we write \linebreak ${\frak A}=\overline{\mathrm{Ass}(a_0,\ldots.a_k)}$

The universal control affine system with $k$-dimensional control is the system
$$
\dot x=x\left(a_0+\sum\limits_{i=1}^ku_ia_i\right), \quad x\in{\frak A},\ u_i\in\R.  \eqno (4)
$$
Given control $u(\cdot)\in L_1([0,t];\R^k)$ and initial condition $x(0)$, we can explicitly write the unique solution of (4)
that is a curve in  ${\frak A}$ whose homogeneous components are absolutely continuous vector functions.

We need some notations in order to write the solution: Let $W$ be the set of words in the alphabet $\{a_0,\ldots,a_k\}$ and
$$
\Delta^n(t)=\{(\tau_1,\ldots,\tau_n\}:0\le \tau_n\le\cdots\le \tau_1\le t\}
$$
be the $n$-dimensional simplex. Given a word $w=a_{i_n}\cdots a_{i_1}$, we set
$$
S^w_u(t)={\int\cdots\int}_{\Delta^n(t)} u_{i_n}(\tau_n)\cdots u_{i_1}(\tau_1)\,d\tau_1\cdots d\tau_n,
$$
where $u_0(t)\equiv 1$. Solutions of (4) have a form:
$$
x(t)=x(0)\sum\limits_{w\in W}S^w_u(t)w.
$$

We keep using chronological notations while working in ${\frak A}$ with the composition ``$\circ$'' substituted by the product in ${\frak A}$. In what follows, we assume that $x(0)=1$. Some notations:
$$
L=\mathrm{Lie}(a_1,\ldots,a_k)\subset \mathrm{Ass}(a_1,\ldots,a_k),\quad {\mathcal L}=\left\{e^V: V\in\overline{L}\right\},
\quad a_u=\sum\limits_{i=1}^ku_ia_i.
$$
Next statement easily follows from basic facts of Lie theory (see \cite{B}), we leave the proof to the reader.

\begin{theorem}
$$
{\mathcal L}=\overline{\left\{\E_0^1a_{v(t)}\,dt: v(\cdot)\in L_1([0,1];\R^k)\right\}}.
$$
Moreover,
$$
\overline{\mathrm{conv}(\mathrm{Ad}{\mathcal L})a_0}=
\overline{\left\{\int_0^1\E_0^t\mathrm{ad}a_{v(\tau)}\,d\tau a_0\,dt : v\in L_1([0,1];\R^k)\right\}}.
$$
$\square$
\end{theorem}

We may translate all the computations made in Section~2 to the universal setting by the substitution of $f_j$ with $a_j$
and we obtain:

\begin{theorem}
For any $t>0$, the attainable set of the system
$$
\dot x=xV, \quad V\in\mathrm{conv}(\mathrm{Ad}{\mathcal L})a_0+L
$$
at $t$ is contained in the closure of the attainable set of system (4) at $t$. \newline
$\square$
\end{theorem}

To go ahead, we need an effective description of $\mathrm{conv}\{(\mathrm{Ad}z)a_0: z\in{\mathcal L}\}$.

Some notations:
$A^n=\bigoplus\limits_{i=0}^nA_n,\ L^n=A^n\cap L,\ \pi^n:{\frak A}\to A^n$
is the continuous projector such that $\pi^n(A_i)=0,\ \forall\,i>n;$ we also set ${\cal L}^n=\pi^n({\cal L})\subset A^n$.

The space $A^n$ is endowed with the truncated products ${x_ny}=\pi^n(xy),\ {[x_ny]}={x_ny}-{y_nx}=\pi^n([x,y])$. Obviously, $\pi^n$ is a homomorphism of ${\frak A}$ on the step $n$ nilpotent associate algebra $A^n$, $\pi^n(xy)={\pi^n(x)}_n\pi^n(y)$. Moreover, $L^n$ is a step $n$ nilpotent Lie algebra and
${\mathcal L}^n$ is the Lie group generated by $L^n$.

Let $V_\alpha,\ \alpha=1,2,\ldots,$ be a linearly ordered homogeneous additive basis of $L$, which satisfies the following property: there exists $r>0$ such that
$$
\min\{n>0: V_\alpha\in L^n\}-\min\{n>0: \pi^n(V_\alpha)=0\}\le r,\quad \forall\,\alpha>0.
$$
Any finite set of linearly independent elements of $L$ is contained in a basis of this type. If $r=1$, then $V_\alpha$ is a homogeneous basis.
It is easy to see that
$$
{\mathcal L}=\left\{\prod\limits_{\alpha=1}^\infty e^{v_\alpha V_\alpha}: v_\alpha\in\R\right\}. \eqno (5)
$$
This is what people call ``the 2nd type coordinates'' for the Lie group, while the presentation ${\mathcal L}=e^{\overline{L}}$ is the ``1st type coordinates''.

We have:
$$
\prod\limits_{\alpha=1}^\infty e^{v_\alpha V_\alpha}=
1+\sum_{\substack{\alpha_1\le\cdots\le\alpha_m\\ i_1,\ldots,i_m>0}}\frac{v_{\alpha_1}^{i_1}\cdots v^{i_m}_{\alpha_m}}{i_1!\cdots i_m!}V^{i_1}_{\alpha_1}\cdots V^{i_m}_{\alpha_m}.
$$
According to the Poincare--Birkhoff--Witt theorem, the elements
$$
V^{i_1}_{\alpha_1}\cdots V^{i_m}_{\alpha_m},\quad \alpha_1\le\cdots\le\alpha_m,\ i_1,\ldots,i_m>0,\ m\ge 0, \eqno (6)
$$
form an additive basis of $\mathrm{Ass}\{a_1,\ldots,a_k\}$.\footnote{The empty product is assumed to be equal to 1.} Next statement is an immediate corollary of this fact.

\begin{prop} $\overline{\mathrm{span}{\mathcal L}}=\overline{\mathrm{Ass}(a_1,\ldots,a_k)}$ and $\mathrm{span}{\mathcal L}^n=A^n,\ \forall\,n\ge 0. \qquad \square$
\end{prop}

Next statement reduces the study of $\mathrm{conv}\{(\mathrm{Ad}z)a_0: z\in{\mathcal L}\}\subset\overline{L}$ to the study
of $\mathrm{conv}{\mathcal L}\subset \overline{\mathrm{Ass}(a_1,\ldots,a_k)}$.

\begin{prop} Linear map $Ad_0:\mathrm{Ass}(a_1,\ldots,a_k)\to\mathrm{Lie}(a_0,a_1,\ldots a_k)$ defined by its action on the basis:
$$
Ad_0\left(V^{i_1}_{\alpha_1}\cdots V^{i_m}_{\alpha_m}\right)=(\mathrm{ad}V_{\alpha_1})^{i_1}\cdots (\mathrm{ad}V_{\alpha_m})^{i_m}a_0
$$
is injective.
\end{prop}
{\bf Proof} Given $v=\{v_\alpha\}_{\alpha=1}^\infty$ we set $Z(v)=\prod\limits_{\alpha=1}^\infty e^{v_\alpha V_\alpha}$.
We have:
$$
Ad_0\left(Z(v)\right)=\prod\limits_{\alpha=1}^\infty e^{v_\alpha\mathrm{ad}V_\alpha}a_0=Z(v)a_0{Z(v)}^{-1}.
$$
Any linear combination of $Ad_0\left(V_{\alpha_1}^{i_1}\cdots V_{\alpha_m}^{i_m}\right)$ can be written as $D_0Ad_0Z(\cdot)$,
where $D_0$ is a linear combination of partial differentials $\frac{\p^{i_1}}{\p v_{\alpha_1}^{i_1}}\cdots\frac{\p^{i_m}}{\p v_{\alpha_m}^{i_m}}\Bigr|_{v=0}$.

Note that $Ad_0\left(Z(v)\right)$ is a series of linear combinations of words where any presented word contains exactly one letter $a_0$. It is easy to see that
$$
\frac{\p^{\bar\alpha}}{\p v^{\bar\alpha}}\Bigr|_{v=0}Ad_0(Z(v))=
\frac{\p^{\bar\alpha}}{\p v^{\bar\alpha}}\Bigr|_{v=0}\left(Z(v)a_0({Z(v)}^{-1}\right)
$$
$$
=\left(\frac{\p^{\bar\alpha}}{\p v^{\bar\alpha}}\Bigr|_{v=0}Z(v)\right)a_0 +\sum\limits_{i=1}^k\xi_ia_i,\quad \xi_i\in{\frak A}.
$$
Hence the equality $D_0Ad_0(Z(\cdot))=0$ implies $\left(D_0Z(\cdot)\right)a_0=0$. It follows that $D_0Z(\cdot)=0$ and
$D_0=0$ because all partial derivatives of $Z(\cdot)$ at 0 are linearly independent.\qquad $\square$

Let $t_\alpha,\ \alpha=1,2,\ldots,$ be coordinates on $L$ induced by the basis $V_\alpha$. In other words, $t_\alpha\in L^*,\
\langle t_\alpha,V_\beta\rangle=\delta_{\alpha,\beta}$. The basis provides the identification of $\bar L$ and $L^*$.
According to this identification, a series $\sum_\alpha v_\alpha V_\alpha$ is identified with the linear function
$\sum_\alpha v_\alpha t_\alpha$ on $L$.

Moreover, the monomials $t^{i_1}_{\alpha_1}\cdots t^{i_m}_{\alpha_m}$ are coordinates on the vector space
$\mathrm{Ass}(a_1,\ldots,a_k)$. A monomial $t^{i_1}_{\alpha_1}\cdots t^{i_m}_{\alpha_m}$ treated as a linear form on \linebreak
$\mathrm{Ass}(a_1,\ldots,a_k)$ annihilates all elements of the basis (6) except of $V^{i_1}_{\alpha_1}\cdots V^{i_m}_{\alpha_m}$ and
$
\langle t^{i_1}_{\alpha_1}\cdots t^{i_m}_{\alpha_m},V^{i_1}_{\alpha_1}\cdots V^{i_m}_{\alpha_m}\rangle=1.
$

The basis (6) provides the identification of $\overline{\mathrm{Ass}(a_1,\ldots,a_k)}$ with \linebreak $\mathrm{Ass}(a_1,\ldots,a_k)^*$
and eventually with the space of formal power series on the variables $t_\alpha,\ \alpha=1,2,\ldots,$. Let ${\mathcal S}$ the space of formal power series and
$$
\nu:\overline{\mathrm{Ass}(a_1,\ldots,a_k)}\to {\mathcal S}
$$
be the continuous isomorphism of vector spaces that realizes the mentioned identification,
$$
\nu:V^{i_1}_{\alpha_1}\cdots V^{i_m}_{\alpha_m}\mapsto t^{i_1}_{\alpha_1}\cdots t^{i_m}_{\alpha_m},
$$
where $\alpha_1\le\cdots\le\alpha_m$ as in (6).
Linear map $\nu$ depends on the choice of the basis $V_\alpha$ and it is not a homomorphism of the algebras.

\begin{definition} We say that a nonzero function $\varphi:L\to\R$ is exponential if the restriction of $\varphi$ to any
finite-dimensional subspace of $L$ is continuous and
$$
\varphi(z_1+z_2)=\varphi(z_1)\varphi(z_2),\quad \forall\,z_1,z_2\in L.
$$
\end{definition}

It is easy to see that, written in the coordinates, exponential functions are exactly functions of the form
$\varphi(t)=e^{\langle v,t\rangle}$, where
$$
v=\{v_\alpha\}_{\alpha=1}^\infty,\quad t=\{t_\alpha\}_{\alpha=1}^\infty,\quad \langle v,t\rangle=\sum\limits_{\alpha=1}^\infty v_\alpha t_\alpha.
$$
Recall that an element of $L$ has only a finite number of nonzero coordinates $t_\alpha$.

The space of exponential functions is denoted by ${\mathcal E}$. The identification of the exponential function with the exponential series gives the inclusion ${\mathcal E}\subset {\mathcal S}$.

\begin{prop} $\nu({\mathcal L})={\mathcal E}$.
\end{prop}

\noindent{\bf Proof.} Indeed,
$$
\nu\left(\prod\limits_{\alpha=1}^\infty e^{v_\alpha V_\alpha}\right)=
1+\sum_{\substack{\alpha_1\le\cdots\le\alpha_m\\ i_1,\ldots,i_m>0}}\frac{v_{\alpha_1}^{i_1}\cdots v^{i_m}_{\alpha_m}}{i_1!\cdots i_m!}\nu\left(V^{i_1}_{\alpha_1}\cdots V^{i_m}_{\alpha_m}\right)=
$$
$$
1+\sum_{\substack{\alpha_1\le\cdots\le\alpha_m\\ i_1,\ldots,i_m>0}}\frac{v_{\alpha_1}^{i_1}\cdots v^{i_m}_{\alpha_m}}{i_1!\cdots i_m!}t^{i_1}_{\alpha_1}\cdots t^{i_m}_{\alpha_m}=
\prod\limits_{\alpha=1}^\infty e^{v_\alpha t_\alpha}=e^{\langle v,t\rangle}. \qquad\square
$$
We see that $\nu({\mathcal L})$ does not depend on the choice of the basis $V_\alpha$, unlikely the isomorphism $\nu$.

We are going to give a special description of $\mathrm{conv}{\mathcal E}(L)$ (and hence of $\mathrm{conv}{\mathcal L}$),
which allows to effectively compute a broad class of finite-dimensional projections of this set. We need the following

\begin{lemma} Let $w_1,\ldots,w_n\in\R^m$ be mutually distinct vectors: $w_i=w_j$ if and only if $i=j$; then the Taylor
polynomials of order $(n-1)$  at 0 of the exponential functions $e^{\langle w_1,t\rangle},\ldots, e^{\langle w_n,t\rangle}$ are linearly independent.
\end{lemma}

\noindent{\bf Proof.} We start from the case $m=1$. In this case, $w_i$ are mutually different scalars and the Taylor polynomials have a form $\sum\limits_{j=0}^{n-1}\frac{w_i^j}{j!}t^j$. The vectors $\left(1,\ldots,\frac{w_i^{n-1}}{(n-1)!}\right)$
are linearly independent because the Vandermonde matrix
$$
\begin{pmatrix}1& 1&\ldots & 1\\ w_1 & w_2 &\ldots & w_n\\ \vdots & \vdots &\ddots & \vdots\\
  w_1^{n-1} & w_2^{n-1} &\ldots & w_n^{n-1} \end{pmatrix}
$$
is nondegenerate.

General case can be easily reduced to the case $m=1$. We simply take
$$
s\in\R^m\setminus\bigcup\limits_{i<j}\{t\in\R^m:\langle w_i,t\rangle=\langle w_j,t\rangle\}
$$
and consider restrictions of the functions $t\mapsto e^{\langle w_i,t\rangle},\ i=1,\ldots,n$, to the line $\R s.\qquad\square$

The space of formal series ${\mathcal S}$ is the adjoint space to the space of (finite) linear combinations of partial differentials $\frac{\p^{i_1}}{\p t_{\alpha_1}^{i_1}}\cdots\frac{\p^{i_m}}{\p t_{\alpha_m}^{i_m}}\bigr|_{t=0}$, where
the pairing of the differential and the series is just the action of the differential on the series.

We set:
$$
{\mathcal D}=\mathrm{span}\left\{\frac{\p^{i_1}}{\p t_{\alpha_1}^{i_1}}\cdots\frac{\p^{i_m}}{\p t_{\alpha_m}^{i_m}}\Bigr|_{t=0}: \alpha_1\le\cdots\le\alpha_m,\ i_1,\ldots,i_m>0,\ m\ge 0\right\},
$$
${\mathcal S}={\mathcal D}^*.$
To any $\varphi\in{\mathcal S}$ we associate a quadratic form $Q_\varphi$ by the following formula:
$$
Q_\varphi(\eta)=\eta_t\eta_s\varphi(t+s),\quad\eta\in{\mathcal D},
$$
where $\eta_t$ differentiates with respect to $t$ and $\eta_s$ differentiates with respect to $s$.

\begin{theorem} Let $\varphi\in\mathrm{span}{\mathcal E},\ \varphi(0)=1$. Quadratic form $Q_\varphi$ is nonnegative
if and only if $\varphi\in\mathrm{conv}{\mathcal E}$.
\end{theorem}

\noindent{\bf Proof.} If $\varphi\in{\mathcal E}$, then $\varphi(t+s)=\varphi(t)\varphi(s)$ and
$$
Q_\varphi(\eta)=\eta_t\eta_s\varphi(t+s)=(\eta \varphi)^2\ge 0.
$$
Let $\varphi=\sum\limits_{i=1}^nc_i\varphi_i$. We may assume that $\varphi_i(t)=e^{\langle w_i,t\rangle},\
i=1,\ldots,n,$ where $w_1,\ldots,w_n$ are mutually distinct. Of course, there exists $m>0$ such that the truncation of these infinite vectors to $\R^m$ are also mutually distinct.

We have: $Q_\varphi(\eta)=\sum\limits_{i=1}^nc_i(\eta\varphi_i)^2$. If all $c_i$ are nonnegative, then $Q_\varphi\ge 0$.
Assume that $c_{i_0}<0$. According to Lemma~1, Taylor polynomials of order $n$ of $\varphi_1,\ldots,\varphi_n$ at 0
are linearly independent. Hence there exists $\eta_0\in{\mathcal D}$ such that $\eta_0\varphi_{i_0}=1,\ \eta_0\varphi_i=0,\ \forall\,i\ne i_0,$ and $Q_\varphi(\eta_0)=c_{i_0}.\qquad \square$

Let $\phi\in\overline{\mathrm{conv}{\mathcal E}}\subset\overline{\mathrm{span}{\mathcal E}}={\mathcal S}$; then
$Q_\phi\ge 0$ by continuity. The inverse statement is not true unfortunately: inequality $Q_\phi\ge 0$ does not imply the relation $\phi\in\overline{\mathrm{conv}{\mathcal E}}$. In other words, the cone $\{\phi\in{\mathcal S}:Q_\phi\ge 0\}$
is strictly bigger than the convex cone generated by ${\mathcal E}$.

To see the difference between these two cones we simply describe the dual cones.
In what follows, we use multi-indices. Given a nonnegative integer $m$, $\mathbb Z^m_+$ is the set of $m$-dimensional vectors
with nonnegative integral coordinates\footnote{$\mathbb Z_+^0=\{0\}$.}. We extend $m$-dimensional vectors by zeros and assume that $\mathbb Z^m_+\subset\mathbb Z^{m'}_+$ if $m\le m'$; then $\mathbb Z_+^\infty=\bigcup\limits_{m\ge 0}\mathbb Z^m_+$ is a set of infinite vectors with a finite number nonzero coordinates. If $i=(i_1,\ldots,i_m)\in\mathbb Z^m_+$, then:
$$
t^i=t_1^{i_1}\cdots t_m^{i_m},\quad |i|=\sum\limits_{j=1}^mi_j,\quad \varphi_0^{(i)}=\frac{\p^{|i|}\varphi}{\p t_1^{i_1}\cdots\p t_m^{i_m}}\Bigr|_{t=0}.
$$

It is convenient to identify operators from ${\mathcal D}$ with their symbols that are linear combinations of monomials
$\tau^i,\ i\in{\mathbb Z}_+^\infty$, where the symbol of $\frac{\p^{|i|}}{\p t^i}$ is $\tau^i$.

The pairing of ${\mathcal D}$ and ${\mathcal S}={\mathcal D}^*$ is defined by the rule:
$\langle\tau^i,t^j\rangle=i!\delta_{ij}$. Given a closed convex cone $K\subset{\mathcal S}$, the dual cone $K^\circ\subset{\mathcal D}$ is defined as follows:
$$
K^\circ=\{\xi\in{\mathcal D}: \xi\varphi\ge 0,\ \,\forall\varphi\in K\}.
$$

\begin{lemma} Dual cone to the closed convex con generated by ${\mathcal E}$ is the set of differential operators whose symbols are nonnegative polynomials.

Dual cone to $\{\phi\in{\mathcal S}:Q_\phi\ge 0\}$ is the set differential operators whose symbols are sums of squares of real polynomials.
\end{lemma}
{\bf Proof.} The statement follows directly from the definitions. Let $\xi=\sum\xi_i\tau^i$ be the symbol of a differential operator from ${\mathcal D}$. The polynomial $\xi$ belongs to ${\mathrm{cone}{\mathcal E}}^\circ$ if and only if
$\sum\xi_iv^i\ge 0,\ \,\forall v$.

Moreover, let $\phi\in{\mathcal S}$, then $Q_\phi\ge 0$ if and only if $\langle\eta^2,\phi\rangle\ge 0$ for any symbol $\eta$. Hence the cone $\{\phi\in{\mathcal S}:Q_\phi\ge 0\}$ and the cone generated by the squares of real polynomial are dual to each other.$\qquad\square$

It is well-known that not any nonnegative polynomial is a sum of squares. By duality, not any $\phi$ with the nonnegative Hankel form $Q_\phi$ belongs to the closed convex cone generated by ${\mathcal E}$.

The space $\overline{\mathrm{span}\{(\mathrm{Ad}z)a_0 :z\in{\mathcal L}\}}$ contains only bracket polynomials of degree 1 with respect to $a_0$. Let us show that all bracket polynomials of degree 1 with respect to $a_0$ are there. Indeed,
the closure of $\mathrm{span}\{(\mathrm{Ad}z)a_0 :z\in{\mathcal L}\}$ is equal to the closure of
$$
\mathrm{span}\{\mathrm{ad}a_{i_1}\cdots(\mathrm{ad}a_{i_n})a_0 :1\le i_1,\ldots,i_n\le k,\ n\ge 0\}.
$$

\begin{prop} Any Lie bracket polynomial in the variables $a_0,\ldots,a_k$ of degree 1 with respect to $a_0$ is a linear combination of the Lie monomials of the form
$$
\mathrm{ad}a_{i_1}\cdots(\mathrm{ad}a_{i_n})a_0,\quad 1\le i_1,\ldots,i_n\le k,\ n\ge 0.
$$
\end{prop}

\noindent{\bf Proof.} As we know, linear map
$$
\pi:\mathrm{Ass}(a_0,\ldots,a_k)\to \mathrm{Lie}(a_0,\ldots,a_k)
$$
defined by its action on the words:
$$
\pi:a_{i_1}\cdots a_{i_n}\mapsto \frac 1n (\mathrm{ad}a_{i_1})\cdots (\mathrm{ad}a_{i_{n-1}})a_{i_n}
$$
is a projection of $\mathrm{Ass}(a_0,\ldots,a_k)$ on $\mathrm{Lie}(a_0,\ldots,a_k)$. In order to prove the proposition, we have to show that the Lie monomial of the form \linebreak $(\mathrm{ad}a_0)(\mathrm{ad}a_{i_1})\cdots (\mathrm{ad}a_{i_{n-1}})a_{i_n}$ is a linear combination of the monomials of the form $(\mathrm{ad}a_{j_1})\cdots(\mathrm{ad}a_{j_n})a_0$.

We use induction with respect to $n$. The case $n=1$ is trivial. Induction step $(n-1)\to n$ is as follows.

$$
(\mathrm{ad}a_0)(\mathrm{ad}a_{i_1})\cdots (\mathrm{ad}a_{i_{n-1}})a_{i_n}=(\mathrm{ad}a_{i_1})(\mathrm{ad}a_0)\cdots (\mathrm{ad}a_{i_{n-1}})a_{i_n}+
$$
$$
+(\mathrm{ad}[a_0,a_{i_1}])\cdots (\mathrm{ad}a_{i_{n-1}})a_{i_n}.
$$
The first term in the right-hand side is fine due to the induction assumption. For the second term, we have:
$$
(\mathrm{ad}[a_0,a_{i_1}])\cdots (\mathrm{ad}a_{i_{n-1}})a_{i_n}=
(\mathrm{ad}a_{i_2})(\mathrm{ad}[a_0,a_{i_1}])\cdots (\mathrm{ad}a_{i_{n-1}})a_{i_n}
$$
$$
+(\mathrm{ad}[[a_0,a_{i_1}],a_{i_2}])\cdots (\mathrm{ad}a_{i_{n-1}})a_{i_n}.
$$
The induction assumption implies that any Lie polynomial of degree $\le n$ and degree 1 with respect to $a_0$ has a desired
presentation. Hence the first term in the right-hand side of the last identity is also fine and we may focus on the term
$
(\mathrm{ad}[[a_0,a_{i_1}],a_{i_2}])\cdots (\mathrm{ad}a_{i_{n-1}})a_{i_n}.
$
We proceed again and again in the same way and eventually reduce everything to the term
$$
(\mathrm{ad}[\cdots[a_0,a_{i_1}],\cdots],a_{i_{n-1}}])a_{i_n}=(-1)^n(\mathrm{ad}a_{i_n})\cdots(\mathrm{ad}a_{i_1})a_0.
\eqno \square
$$

Theorem~2 implies that any element of the set $\overline{\mathrm{conv}(Ad{\mathcal L})a_0}+\overline{L}$ provides a good bracket. It happens that this set is actually the set of all good brackets for the
universal system (4).

\begin{theorem} Let $V\in\overline{Lie(a_0,\ldots,a_k)}$; the left-invariant vector field $xV$ is a good bracket for system (4) if and only if
$V\in \overline{\mathrm{conv}(\mathrm{Ad}{\mathcal L})a_0}+\overline{L}$
\end{theorem}

\noindent{\bf Proof.} If $V\in \overline{\mathrm{conv}(\mathrm{Ad}{\mathcal L})a_0}+\overline{L}$ then, according to Theorem~2, $xV$ is indeed a good bracket for system (4).

Now let $V\notin \overline{\mathrm{conv}(\mathrm{Ad}{\mathcal L})a_0}+\overline{L}$. We are going to show that $e^{tV}$ does not belong to the closure of the attainable sets of system (4) for sufficiently small $t>0$. First of all, there exists $n>0$
such that
$$
\pi^n(V)\notin \pi^n\left(\overline{\mathrm{conv}(\mathrm{Ad}{\mathcal L})a_0}+\overline{L}\right)=
\overline{\mathrm{conv}(Ad{\pi^n(\mathcal L}))a_0}+\pi^n\overline{L})
$$
and we can work with the truncated system $(\pi^nx)^\cdot=(\pi^nx)\pi^nV$ on the finite-dimensional Lie group $G\subset A^n$, the exponential of the
nilpotent Lie algebra $\pi^n\mathrm{Lie}(a_0,\ldots,a_k)$. In what follows, until the end of the section, we keep notations $x,V,{\frak A},{\mathcal L}$ etc.
for $\pi^nx,\pi^nVV,\pi^n{\frak A},\pi^n{\mathcal L},\ldots$ in order to simplify the presentation.
Any big enough $n$ is fine and we simply omit symbol $\pi^n$.

Let $K_0=\overline{\mathrm{conv}(\mathrm{Ad}{\mathcal L})a_0}$ and $K=\mathbb R_+K_0$, the closed convex cone generated by $(\mathrm{Ad}{\mathcal L})a_0$. Let $H_0$ be affine hyperplane in $A^n$ which consists of the linear combinations of
words such that the coefficient against $a_0$ equals 1. Obviously, $K_0=K\cap H_0$.

The cone $K$ is acute, it does not contain any proper subspace; indeed, Lemma~2 implies the dual cone $K^\circ$ has a nonempty interior. We consider the system
$$
\dot x=xu,\quad x\in G,\ u\in K+L,  \eqno (7)
$$
and we are going to show that $e^{tV}$ does not belong to the closure of the attainable set of system (7) for sufficiently small $t>0$. Of course, it is enough to consider solutions of system (7) with the initial condition x(0)=1; all other solutions are left translation of these ones.

The cone $K$ is invariant with respect to the adjoint action of ${\mathcal L}$,
$\mathrm{Ad}{\mathcal L}(K)=K$, as well as the space $L$. Hence system (7) has a well-defined
projection on the homogeneous space $G/{\mathcal L}=\{x{\mathcal L}: x\in G\}.$

Let us describe the projected system. We have :
$$
\mathrm{Ass}(a_0,\ldots,a_k)=\mathrm{Ass}(a_1,\ldots,a_k)\oplus{\mathcal A}^0,
$$
where ${\mathcal A}^0$ is the linear hull of words which contain the letter $a_0$. Similarly, $\mathrm{Lie}(a_0,\ldots,a_k)=L\oplus L^0$, where $L^0\subset{\mathcal A}^0$; then $L^0$ is an ideal of $\mathrm{Lie}(a_0,\ldots,a_k)$ and
${\mathcal L}^0\doteq e^{L^0}$ is a normal subgroup of $G$. Moreover, the map
$$
y\mapsto y{\mathcal L},\quad y\in{\mathcal L}^0
$$
is a diffeomorphism of ${\mathcal L}^0$ on the $G/{\mathcal L}$. The projected system takes the form:
$$
\dot y=yv, \quad y\in {\mathcal L}^0,\ v\in K,\ y(0)=1. \eqno (8)
$$

Attainable sets of system (8) are the projections of the attainable sets of system (7). Let $\gamma(t)\in{\mathcal L}^0$
be the projection of $e^{tV}$, then$\dot\gamma\notin K$. We have to show that $\gamma(t)$ does not belong to the attainable sets of system (8) for small positive $t$.

\begin{lemma} Let $K$ be a closed convex acute cone in $\mathbb R^m,\ z\in\mathbb R^m\setminus K$, and $B(q)$ be a smooth family of $m\times m$ matrices, $B(0)=Id$. Then there exists $c>0$ and a neighborhood of the origin $O\subset\mathbb R^m$
such that the relations
$$
q(0)=0,\quad \dot q(\tau)\in B(q(\tau))K,\quad q(\tau)\in O,\quad 0\le\tau\le t,
$$
imply the inequality $|q(t)-tz|\ge ct$.
\end{lemma}

\noindent{\bf Proof.} There exists $\xi\in\mathrm{int}K^\circ\subset{\mathbb R^m}^*$ such that $\langle\xi,z\rangle<0$.
Then $\langle\xi,y\rangle>0$ for any $y\in K\setminus\{0\}$ and, moreover, there exists a neighborhood of the origin $O\subset\mathbb R^m$ such that $\langle\xi,B(q)y\rangle>0$
for any $y\in K\setminus\{0\}$ and any $q\in O$. Hence $\langle\xi,\dot q(\tau)\rangle\ge 0,\ 0\le\tau\le t,$ and
$\langle\xi,q(t)\rangle\ge 0,\qquad\square$

Let us take local coordinates in a neighborhood of the unit element in ${\mathcal L}^0$ such that the coordinate presentation of $\gamma(\cdot)$ is a linear vector function, $\gamma(t)=t\dot\gamma(0)$ and apply the lemma. We obtain that an admissible trajectory of system (9) $q(t)$ is separated from $\gamma(t)$ as long as $q(\tau)\in O,\ \forall\,\tau\in[0,t]$.

Unfortunately, this local statement is not sufficient for our purposes: the cone $K$ is unbounded and an admissible trajectory may go far from the starting point in short time and then, who knows, to come back (from the other side of the
hyperplane $\xi^\bot$). It may happen for general system. We are going to show that this never happens for system (8) and
thus to complete the proof of the theorem.

We assume that $\mathbb R\oplus{\mathcal A}^0$ is equipped with an Euclidean structure. Given
$y\in{\mathcal L}^0\subset\mathbb R\oplus{\mathcal A}^0$ we denote by $|y|$ the length of $y$ in
$\mathbb R\oplus{\mathcal A}^0$.

\begin{lemma} There exists $c>0$ such that $|y(t)-1|\ge c\int\limits_0^t|v(\tau)|\,d\tau$ for any integrable
$v(\cdot):[0,t]\mapsto K$, where
$
\dot y(\tau)=y(\tau)v(\tau),\ 0\le\tau\le t,\ y(0)=1.
$
\end{lemma}

\noindent{\bf Proof.} Let $\xi\in\mathrm{int}K^\circ\subset{{\mathcal A}^0}^*$, then there exists $\delta>0$ such that
$\langle\xi,v\rangle\ge\delta|v|$ for any $v\in K$.

We have: ${\mathcal A}^0=\bigoplus\limits_{i\ge 1}A^0_i$, where $A^0_i$ is the linear hull of all words of degree $i$ with respect to $a_0$. The cone $K$ is contained in $A^0_1$; hence we can take $\xi\in\mathrm{int}K^\circ$ such that
$\xi\perp\bigoplus\limits_{i>1}A^0_i$. Moreover,
$$
y(t)=1+\int\limits_0^ty(\tau)v(\tau)\,d\tau=1+\int\limits_0^tv(\tau)\,d\tau+Y(t),
$$
where $Y(t)\in \bigoplus\limits_{i>1}A^0_i$. Hence
$$
\langle\xi,y(t)-1\rangle=\int\limits_0^t\langle\xi,v(\tau)\rangle\,d\tau\ge\delta\int\limits_0^t|v(\tau)|\,d\tau
\eqno \square$$

We are now ready to complete the proof of the theorem. Indeed, there exists $\e>0$ such that the inequality
$\int\limits_0^t|v(\tau)|\,d\tau\le\e$ implies that $y(\tau),\ 0\le\tau\le t,$ belongs to the neighborhood $O$ guaranteed by Lemma~3 and $y(t)$ is separated from $e^{tV}$. On the other hand, the inequality $\int\limits_0^t|v(\tau)|\,d\tau\ge\e$
implies that $|y(t)-1|\ge\e c$ and $y(t)$ is again separated from $e^{tV}$ for small $t.\qquad \square$

\section{Truncation}

Free Lie algebra is too big. It is more practical and sufficient for many purposes to consider its finite-dimensional nilpotent truncations.

Let $I\subset\mathbb Z_+^\infty$ be a finite subset, we set:
$$
{\mathcal P}(I)=\left\{\sum\limits_{i\in I}c_it^i:c_i\in\R\right\},
$$
a $\#I$-dimensional space of polynomials. We denote by $\Pi_I:{\mathcal S}\to{\mathcal P}(I)$ the continuous linear projector defined by the rule:
$$
\Pi_I(t^i)=\begin{cases}t^i,&\text{if $i\in I$;}\\ 0,&\text{if $i\in\mathbb Z_+^\infty\setminus I$.} \end{cases}
$$

\begin{lemma}$\mathrm{span}\Pi_I({\mathcal E})={\mathcal P}(I)$, for any $I\subset\mathbb Z_+^\infty$ such that $\#I<\infty$.
\end{lemma}

\noindent{\bf Proof.} Given positive integers $m,n$, we consider the ``box'':
$$
J^m_n=\{(i_1,\ldots,i_m): 0\le i_j<n,\ j=1,\ldots,m\}\subset\mathbb Z_+^m.
$$
Any finite set is contained in some box, $I\subset J^m_n$. Let $v_,\ldots,v_n$ be mutually distinct real numbers.
I claim that $n^m$ vectors
$$
\Pi_{J^m_n}\Bigl(e^{\sum\limits_{j=1}^mv_{i_j}t_j}\Bigr),\quad 1\le i_j\le n,\ j=1,\ldots,m,
$$
form a basis of ${\mathcal P}(J^m_n)$.

Indeed, ${\mathcal P}(J^m_n)={\mathcal P}(J^1_n)^{\otimes m}$. According to Lemma~1, the vectors
$\Pi_{J^1_n}\left(e^{v_it}\right),\ i=1,\ldots,n,$ form a basis of ${\mathcal P}(J^1_n)$. Moreover,
$$
\Pi_{J^m_n}\Bigl(e^{\sum\limits_{j=1}^mv_{i_j}t_j}\Bigr)=\Pi_{J^1_n}\left(e^{v_{i_1}t}\right)\Bigr|_{t=t_1}\cdots
\Pi_{J^1_n}\left(e^{v_{i_m}t}\right)\Bigr|_{t=t_m}
$$
is just tensor product of $\Pi_{J^1_n}\left(e^{v_{i_j}t}\right),\ j=1,\ldots,m$ and tensor products of the elements of the basis form the basis of the tensor product of the spaces.

It follows that $\mathrm{span}\Pi_{J^m_n}({\mathcal E})={\mathcal P}(J^m_n)$; hence
$\mathrm{span}\Pi_I({\mathcal E})={\mathcal P}(I)$, for any $I\subset J^m_n.\qquad \square$

The next result gives a useful estimate of $\Pi_I({\mathrm{conv}\mathcal E})$ for a broad class of finite subsets $I$.
Let $C\subset\R^m_+$ be a convex compact subset; we set: $I_C=C\cap\mathbb Z_+^m$.

\begin{theorem} Let $C$ be a convex compact subset of $\R^m_+$ and $0\in C$. If $\phi\in{\mathcal P}(I_C)$ belongs to $\overline{\Pi_{I_C}({\mathrm{conv}\mathcal E})}$, then
$$
\phi(0)=1,\quad \sum\limits_{i,j\in I_{\frac{1}{2}C}}\phi_0^{(i+j)}\xi_i\xi_j\ge 0,\quad \forall\,\xi_i\in\R,\, i\in I_{\frac{1}{2}C}. \eqno (9)
$$
Moreover, if $m=1$ or $|i|\le 2,\,\forall i\in I_C$, then condition (9) is sufficient for $\phi$ to belong to
$\overline{\Pi_{I_C}({\mathrm{conv}\mathcal E})}$.
\end{theorem}

\noindent{\bf Proof.} Let $I\subset\mathbb Z^m_+,\ \#I<\infty$; we set:
$
{\mathcal D}_I=\mathrm{span}\left\{\frac{\p^{|i|}}{\p t^i}\Bigr|_{t=0}\right\}\subset{\mathcal D}.
$

Let $\phi=\Pi_{I_C}(\varphi),\ \varphi\in\mathrm{conv}{\mathcal E}$. The restriction of the form $Q_\varphi$ to the subspace
${\mathcal D}_{I_{\frac 12C}}$ depends only on $\phi$; indeed, if $i,j\in\frac 12C$, then $2i,2j\in C$, and $i+j\in C$
due to the convexity of $C$. Quadratic form (9) is equal to $Q_\varphi\bigr|_{{\mathcal D}_{I_{\frac 12C}}}$.

The form $Q_\varphi$ is nonnegative (see Theorem~3), hence its restriction is nonnegative and we obtain that the nonnegativity of the quadratic form (6) is a necessary condition for $\phi$ to belong to $\Pi_{I_C}({\mathcal E})$.

In general, this condition is not sufficient: the cone $\bigl\{\phi\in{\mathcal P}(I_C): Q_\phi\bigr|_{{\mathcal D}_{I_{\frac 12C}}}\ge 0\bigr\}$ may be strictly bigger than the closed convex cone generated by $\Pi_{I_C}({\mathcal E})$. These two cones are equal if the dual cones are equal, i.e. if any nonnegative polynomial from ${\mathcal P}(I_C)$ is a sum of squares of polynomials from
${\mathcal P}\bigl(I_{\frac 12C}\bigr)$ (see Lemma~2). The last property is valid for polynomials of one variable and for quadratic polynomials.$\qquad \square$

\medskip Let us show that any $\phi\in{\mathcal P}_{I_C}$ such that $Q_\phi\bigr|_{{\mathcal D}_{I_{\frac 12C}}}>0$ can be extended to $\varphi\in{\mathcal S}$ such that $Q_\varphi>0$.
Let $B_r\subset\R^m$ be the ball of radius $r$ centered at 0. We set: $C^r=\frac 12C+B_r$, the radius $r$ neighborhood of $\frac 12C$.

\begin{lemma} Let $\phi\in{\mathcal P}(I_C)$ be such that $Q_\phi\bigr|_{{\mathcal D}_{I_{\frac 12C}}}>0$; then for any $r\ge 0$ there exists $\phi_r\in{\mathcal P}(I_{2C^r})$ such that $\Pi_{I_C}(\phi_r)=\phi$ and $Q_{\phi_r}\bigr|_{{\mathcal D}_{I_{C^r}}}>0$. Moreover, there exists $\varphi\in{\mathcal S}$ such that $\Pi_{I_C}(\varphi)=\phi$ and $Q_\phi>0$.
\end{lemma}

\noindent{\bf Proof.} We construct $\phi_r$ step by step starting from $\phi_0=\phi$. Assume that $\phi_r$ is already constructed for some $r$. Let
$$
r'=\min\{\rho : \exists\,i\in 2C^r\setminus C^r : \text{such that}\ i\in C^\rho\}.
$$
Let $J=\{i\in 2C^r\setminus C^r : i\in C^{r'}\}$, $\eta>0$. We define $\phi^\eta_r\in{\mathcal P}\left(I_{2C^{r'}}\right)$ by the formula:
$$
\phi^\eta_r(t)=\phi_r(t)+\eta\sum\limits_{i\in J}\frac{t^i}{i!}.
$$
Then, due to our choice of $r'$ and strict convexity of $C^{r'}$ we obtain:
$$
Q_{\phi^\eta_r}\bigr|_{{\mathcal D}_{I_{C^{r'}}}}(\xi)=\sum\limits_{i+j\in I_{2C^r}}\phi_{r0}^{(i+j)}\xi_i\xi_j+\eta\sum\limits_{i\in J}\xi_i^2. \eqno (10)
$$

I claim that form (10) is positive definite for $\eta$ big enough. Indeed, assume that, on the contrary, for any $N>0$
there exists $\xi^N\in{\mathcal D}_{I_{C^{r'}}},\  \xi^N\ne 0$, such that $Q_{\phi^N_r}(\xi^N)\le 0$. We may assume that
$|\xi^N|=1$. Let $\bar\xi$ be a limiting point of the family $\xi^N$ as $N\to\infty$. We have:
$$
0\ge\frac 1N Q_{\phi^N_r}(\xi^N)=\frac 1N\sum\limits_{i+j\in I_{2C^r}}\phi_{r0}^{(i+j)}\xi^N_i\xi^N_j+\sum\limits_{i\in J}(\xi^N_i)^2,\quad \forall\,N.
$$

It follows that $\bar\xi_i=0,\ \forall\,i\in J$. Hence $\bar\xi\in{\mathcal D}_{I_{C^r}}$ and $Q_{\phi_r}(\bar\xi)>0$.
On the other hand,
$$
0\ge Q_{\phi^N_r}(\xi^N)\ge \sum\limits_{i+j\in I_{2C^r}}\phi_{r0}^{(i+j)}\xi^N_i\xi^N_j,\quad \forall\,N.
$$
Hence $Q_{\phi_r}(\bar\xi)\le 0$.
This contradiction justifies the claim and we take \linebreak $\phi_{r'}=\phi_r^\eta$ for some $\eta$ big enough.

Now, going step by step we reach any value of the radius in a finite number of steps.
What remains is to take $\varphi=\lim\limits_{r\to\infty}\phi_r. \qquad\square$

\medskip
Given $I\subset\mathbb Z^m_+,\ \#I<\infty$, and $\phi\in{\mathcal P}(I)$, we set
$$
erank(\phi)=\min\Bigl\{k>0: \phi=\Pi_I\Bigl(\sum\limits_{i=1}^kc_ie^{\langle\nu_i,t\rangle}\Bigr),\ \nu_i\in\mathbb R^m, c_i\in\mathbb R \Bigr\}.
$$

Let $T_n^m=\{i\in \mathbb Z^m_+: |i|<n\}$. The following dichotomy helps to distinguish projections of $\mathrm{conv}{\mathcal E}$ from other polynomials $\phi\in{\mathcal P}(I_C)$ such that
$Q_\phi\bigr|_{{\mathcal D}_{I_{\frac 12C}}}\ge 0$.

\begin{prop} Let $\phi\in{\mathcal P}(I_C),\ Q_\phi\bigr|_{{\mathcal D}_{I_{\frac 12C}}}\ge 0$, and $I_C\subset T^m_n$;
\begin{enumerate}
\item if $\phi\notin\mathrm{conv}{\mathcal E}$ and $\phi=\Pi_{I_C}(\phi_n)$, where $\phi_n\in{\mathcal P}(T^m_{2n}),\
Q_\phi\bigr|_{{\mathcal D}_{T^m_n}}\ge 0$, then $erank(\phi_n)>n$;
\item if $\phi\in\mathrm{conv}{\mathcal E}$, then there exists $\phi_n\in{\mathcal P}(T^m_{2n})$ such that
$\phi=\Pi_{I_C}(\phi_n),\ Q_\phi\bigr|_{{\mathcal D}_{T^m_n}}\ge 0$ and $erank(\phi_n)\le\#I_C$.
\end{enumerate}
\end{prop}

\noindent{\bf Proof.} Statement 1 follows from Lemma~1 and the proof of Theorem~3. Under the assumption of statement~2,
there exists
$$
\varphi=\sum\limits_{i=1}^{\#I_C}\alpha_ie^{\langle\nu_i,t\rangle},\quad \nu_i\in\mathbb R^m,\ \alpha_i\ge 0,\ \sum\alpha_i=1,
$$
such that $\phi=\Pi_{I_C}(\varphi)$. Indeed, any element of the convex hull of a connected subset in $\mathbb R^k$ is a convex combination of $k$ elements of the subset. We set $\phi_n=\Pi_{T^m_n}(\varphi). \qquad\square$

\section{Iteration}

Good brackets of the universal system are those Lie bracket series that are good for {\it any} system for which they converge and, in particular, for any nilpotent system.

If Lie algebra of the system is not universal, then it may admit more good brackets. Any Lie algebra with $k+1$
generators is the quotient algebra of $\mathrm{Lie}(a_0,\ldots,a_k)$ with respect to some Lie ideal. Here we consider ideals
with a finite number of generators.

Let ${\mathfrak I}\subset \mathrm{Lie}(a_0,\ldots,a_k)$ be such an ideal, then $\mathfrak{AI}=\mathfrak{IA}$ is an ideal of the associate algebra ${\mathfrak A}$. We set:
$$
\hat{\mathfrak A}={\mathfrak A}/\mathfrak{AI},\ \hat G=G/{\mathfrak I},\ \hat L=L+{\mathfrak I}\subset\hat{\mathfrak A},\
\hat K=K+{\mathfrak I}\subset\hat{\mathfrak A},\ \hat{\mathcal L}=e^{\hat L},\ \mathrm{etc.}
$$

As we know, the elements of $K_0=\overline{\mathrm{conv}(\mathrm{Ad}{\mathcal L})a_0}$ provide good brackets for the universal
system (4). Hence any element of $\hat K_0=\overline{\mathrm{conv}(\mathrm{Ad}\hat{\mathcal L})\hat a_0}$ provides a good
bracket for the projection of the universal system to $\hat {\mathfrak A}$. The projected system extended by these good brackets
has a form
$$
\dot x=xu,\quad x\in\hat{\mathfrak A},\quad u\in\hat K_0+\hat L. \eqno (11)
$$

I recall that the cone $K=\mathbb R_+K_0$ is acute and $K_0$ does not contain nontrivial affine subspaces. However, the projected set $\hat K_0$ may, of course, contain a nontrivial affine subspace. Assume that
$$
\hat a_0+\mathrm{span}\{b_1,\ldots,b_l\}\subset\hat K_0
$$
and ${\mathcal B}$ is the subgroup of $\hat G$ generated by $e^{sb_1},\ldots,e^{sb_l},\ s\in\mathbb R$. Then any element of
$\overline{\mathrm{conv}(\mathrm{Ad}{\mathcal B})\hat a_0}$ provides a good bracket for system (11) and hence for the projection of the universal system (4) to $\hat{\mathfrak A}$.

Moreover, if $\overline{\mathrm{conv}((\mathrm{Ad}{\mathcal B})\hat a_0\cup\hat K_0)}$ contains a bigger affine subspace,
then we can iterate the procedure and find more good brackets for the projection of the universal system to $\hat{\mathfrak A}$, etc.

We illustrate the construction with a useful in applications class of ideals with one generator, which gives huge spaces of
good brackets. Let $V\in L,\ m>0$. We consider the ideal $\mathfrak I\subset\mathrm{Lie}(a_0,\ldots,a_k)$ generated by
$(\mathrm{ad}V)^{2m}a_0$,
$$
\mathfrak I=\mathbb R(\mathrm{ad}V)^{2m}a_0+\left[(\mathrm{ad}V)^{2m}a_0,\mathrm{Lie}(a_0,\ldots,a_k)\right].  \eqno (12)
$$
Lie algebra $\mathrm{Lie}(a_0,\ldots,a_k)/\mathfrak I$ is maximal Lie algebra with $k+1$ generators such that the bracket
polynomial $(\mathrm{ad}V)^{2m}a_0$ vanishes.

\begin{prop} Assume that the ideal $\mathfrak I$ is defined by the formula (12) and $\hat V$ is the projection of $V$ to
$\hat L$. Then
$$
\hat a_0+\mathrm{span}\{(\mathrm{Ad}\hat{\mathcal L})(\mathrm{ad}\hat V)^{2m-1}\hat a_0\}\subset \hat K_0.
$$
\end{prop}

\noindent{\bf Proof.} According to the definition of $\mathfrak I$, we have:
$e^{s\hat V}=\sum\limits_{i=0}^{2m-1}\frac{s^i}{i!}\hat V^i$. Let $z\in\hat{\mathcal L}$, then
$(\mathrm{Ad}ze^{\tau^{\frac 1{2m-1}} s\hat V})\hat a_0\subset\hat K_0$ for any $s\in\mathbb R,\,\tau>0$. Moreover,
$$
(\mathrm{Ad}ze^{\tau^{\frac 1{2m-1}} s\hat V})\hat a_0=(\mathrm{Ad}z)e^{\tau^{\frac 1{2m-1}} s\mathrm{ad}\hat V}\hat a_0=
\sum\limits_{i=0}^{2m-1}\frac{\tau^{\frac i{2m-1}}s^i}{i!}(\mathrm{Ad}z)(\mathrm{ad}\hat v)^i\hat a_0.
$$
The convex combination
$$
\left(1-\frac 1{|s|^{2m-1}}\right)\hat a_0+\frac 1{|s|^{\frac 1{2m-1}}}(\mathrm{Ad}ze^{\tau^{\frac 1{2m-1}} s\hat V})\hat a_0
$$
also belongs to $\hat K_0$.
We take limits of this convex combination as $s\to\pm\infty$ and obtain that
$\hat a_0\pm\tau(\mathrm{Ad}z)(\mathrm{ad}\hat V)^{2m-1}\hat a_0$ belongs to $\hat K_0$. In other words, the affine line
$\hat a_0+\mathbb R(\mathrm{Ad}z)(\mathrm{ad}\hat V)^{2m-1}\hat a_0$ is contained in $\hat K_0$ for any
$z\in\hat{\mathcal L}$. Hence
$$
\hat a_0+\mathrm{span}\{(\mathrm{Ad}\hat{\mathcal L})(\mathrm{ad}\hat V)^{2m-1}\hat a_0\}\subset \hat K_0. \eqno \square
$$

\section{Examples}

Here we apply abstract results to some particular classes of systems. Each time we write an original control affine system and a compatible extended system. We say that an extended system is compatible with the original one if all vector fields in
the right-hand side of the extended system are good brackets for the original system.

\subsection{Scalar control}

We consider a control affine system of the form:
$$
\left\{\begin{array}{rlc}\dot x&=&u\\ \dot y&=&\psi(x,y),\\ \end{array}\right. \quad u,x\in\R,\ y\in\R^n,  \eqno (13)
$$
where $\psi$ is a degree $m$ vector-polynomial w.r.t. $x$.

The following system with $(m+1)$-dimensional control is compatible with system (13):
$$
\left\{\begin{array}{rlc}\dot x&=&u\\ \dot y&=&\sum\limits_{i=0}^n\frac{u_i}{i!}\frac{\p^i\psi}{\p x^i}(x,y),\\ \end{array}\right.
\quad
u_0=1,\ \sum\limits_{0\le i,j\le \frac m2}u_{i+j}\xi_i\xi_j\ge 0,\ \forall\,\xi_\cdot \ .
$$

Indeed, let $f_0=\langle\psi,\frac{\p}{\p y}\rangle,\ f_1=\frac{\p}{\p x},\ q=(x,y)$, then (13) has a form $\dot q=f_0+uf_1$.
Moreover, $e^{s\mathrm{ad}f_1}f_0(q)=\langle\psi(x+s,y),\frac{\p}{\p y}\rangle$. It remains to apply Theorems 2 and 5.

\subsection{Step 3 nilpotent algebra}

We consider system (1) and assume that vector fields $f_0,\ldots,f_k$ generate a step~3 nilpotent Lie algebra.
In this case, the following system with \linebreak $\frac 16k(k+1)(2k+7)$-dimensional control is compatible with (1):
$$
\dot q=f_0+\sum\limits_{i=1}^k\left(u_if_i+u_{i0}[f_i,f_0]+\frac{u_{ii}}2[f_i,[f_i,f_0]]\right)+
$$
$$
\sum\limits_{1\le i<j\le k}\left(u_{ij}[f_i,[f_j,f_0]]+v_{ij}[[f_i,f_j],f_0]+w_{ij}[f_i,f_j]\right)+
\sum_{\substack{1\le i<j\le k\\i\le\iota\le k}}w_{ij\iota}[[f_i,f_j],f_\iota],
$$
where
$
\sum\limits_{i,j=0}^ku_{ij}\xi_i\xi_j\ge 0,\ \forall\,\xi_\cdot,\ u_{00}=1,\ u_{ij}=u_{ji},
$
and $v_{ij},w_{ij},w_{ij\iota}\in\R$ are free.

This is one more application of Theorems 2 and 5. Indeed, linear hull of the vectors
$
f_i,\ [f_i,f_j],\ [[f_i,f_j],f_\iota],\ 1\le i<j\le k,\ i\le\iota\le k,
$
is equal to $L=\mathrm{Lie}(f_1,\ldots,f_k)$, while the terms of the sum which involve $f_0$ represent $\mathrm{conv}(\mathrm{Ad}e^L)f_0$, according to Theorem~5.

\subsection{Generalized Kalman test}

We consider a system of the form:
$$
\dot x=f(x)+u, \quad x\in\R^n,\ u\in U, \eqno (14)
$$
where $f$ is a vector polynomial of degree $2m-1$ and $U$ is a vector subspace of $\R^n$.

We define a map $\Phi:\R^n\to\R^n$ and subspaces $V_j\subset\R^n$ as follows:
$$
\Phi(v)=\frac{d^{2m-1}}{d\tau^{2m-1}}f(\tau v)\bigr|_{\tau=0},\quad V_1=U,\ V_{j+1}=\mathrm{span}\Phi(V_j)+V_j.
$$

\begin{prop}
System $\dot x=f(x)+v,\ v\in V_n,$ is compatible with system (14).
\end{prop}

\noindent Classical Kalman test corresponds to $m=1$.

\noindent{\bf Proof.} We can write system (14) in the standard form where $f$ is a coordinate presentation of the field $f_0$,
$k=\dim U$, and $f_1,\ldots,f_k$ are constant vector fields, a basis of the subspace $U$. Given a constant vector field $v$
and $j>0$ we have:
$$
(\mathrm{ad}v)^jf(x)=\frac{d^j}{d\tau^j}f(x+\tau u)\bigr|_{\tau=0}.
$$
Moreover, $(\mathrm{ad}v)^{2m-1}f$ is a constant field and $(\mathrm{ad}v)^{2m}f=0$.

Proposition~8 implies that, for any $u_0\in U$, the system
$$
\dot x=f(x)+v,\quad v\in\, U+\R\Phi(u_0),
$$
is compatible with system (14).
A multiple repetition of the above extension completes the proof.$\qquad\square$.

\end{document}